\documentclass[12pt,a4paper]{article}

\usepackage{amssymb}

\usepackage{amsmath}

\begin{document}
\title{Gr\"{o}bner-Shirshov bases and PBW theorems\footnote{Supported by the NNSF of China (11171118), the
Research Fund for the Doctoral Program of Higher Education of China
(20114407110007), the NSF of Guangdong Province (S2011010003374) and
the Program on International Cooperation and Innovation, Department
of Education, Guangdong Province (2012gjhz0007).}}
\author{
L. A. Bokut\footnote {Supported by RFBR 12-01-00329,
LSS--3669.2010.1 and SB RAS Integration grant No.2009.97 (Russia)
and Federal Target Grant ¡°Scientific and educational personnel of
innovation Russia¡± for 2009-2013
(government contract No.02.740.11.5191).} \\
{\small \ School of Mathematical Sciences, South China Normal
University}\\
{\small Guangzhou 510631, P. R. China}\\
{\small Sobolev Institute of Mathematics, Novosibirsk 630090, Russia}\\
{\small  bokut@math.nsc.ru}\\
\\
 Yuqun
Chen\\
{\small \ School of Mathematical Sciences, South China Normal
University}\\
{\small Guangzhou 510631, P. R. China}\\
{\small yqchen@scnu.edu.cn}}

\date{}

\maketitle \noindent\textbf{Abstract:} We review some applications
of  Gr\"{o}bner-Shirshov bases, including PBW theorems, linear bases
of free universal algebras, normal forms for groups and semigroups,
extensions of groups and algebras, embedding of algebras.



\section{Introduction}

Combinatorial algebra seminar at South China Normal University was
organized by the authors in March 2006. Since then, there were some
30 Master Theses and 4 PhD Theses, about 40 published papers in JA,
IJAC, PAMS, JPAA, Comm. Algebra, Algebra Coll., Siberian Math. J.,
Science in China and other Journals and Proceedings. There were
organized 2 International Conferences (2007, 2009) with E. Zelmanov
as Chairman of the Program Committee and several Workshops. We are
going to review some of the papers.

Our main topic is Gr\"{o}bner-Shirshov bases method for different
varieties (categories) of linear ($\Omega$-) algebras over a field
$k$ or a commutative algebra $K$ over $k$: associative algebras
(including group (semigroup) algebras), Lie algebras, dialgebras,
conformal algebras, pre-Lie (Vinberg right (left) symmetric)
algebras, Rota-Baxter algebras, metabelian Lie algebras,
$L$-algebras, semiring algebras, category algebras, etc. There are
some applications particularly to new proofs of some known theorems.

 As it is well known, Gr\"{o}bner-Shirshov (GS for short) bases
method for a class of algebras based on a Composition-Diamond lemma
 for the class. A general form of a Composition-Diamond Lemma
over a field $k$ is as follows.

\vspace{3mm}

 \noindent{\bf
Composition-Diamond lemma for a class of algebras} Let $M(X)$ be a
free algebra of a category $M$ of algebras over a field $k$, $(N(X),
\leq)$ a linear basis (normal words) of $M(X)$ with an ``admissible"
well order and $S\subset M(X)$. Let $Id(S)$ be the ideal of $M(X)$
generated by $S$ and $\bar s$ the leading term of the polynomial
$s$. Then the following statements are equivalent.

\vspace{3mm}

(i) $S$ is a GS basis (i.e., each ``composition" of polynomials from
$S$ is ``trivial").

\vspace{3mm}

(ii) If $f\in Id(S)$, then the maximal word of $f$ has a form $\bar
f = (a\bar sb),\ s\in S,\ a,b\in N(X)$.

\vspace{3mm}

(iii) $Irr(S) = \{u\in N(X)| u\neq (a \bar sb), s\in S, a,b\in
N(X)\}$ is a linear basis of $M(X|S) = M(X)/Id(S)$.

\ \

There are two kinds of compositions, a composition $(s_1,s_2)_w$ of
two polynomials relative to $w=``lcm"(\bar s_1, \bar s_2)\in N(X)$
and a composition $Com_w(s)$ of one polynomial relative to $w=
\overline{(vs)}$ or $w=\overline{(sv)},\ v,\ w\in N(X)$. Namely, for
monic polynomials $s_1,s_2$,
$$
(s_1,s_2)_w=  lcm(\bar s_1,\bar s_2)|_{\bar s_1\rightarrow s_1} -
lcm(\bar s_1,\bar s_2)|_{\bar s_2\rightarrow s_2},
$$
and $Com_w(s) = vs$ or $sv$ correspondingly.

For example, for words $\overline{s_1}=u,\ \overline{s_2}=v\in X^*,\
lcm(u,v) \in \{ucv,\ c\in X^* \mbox{ (a trivial lcm)};\ u=avb,\
a,b\in X^* \mbox{ (an inclusion lcm)}; ub=av,\ a,b\in X^*,\
deg(ub)<deg(u)+deg(v) \mbox{ (an intersection lcm)}\}$. In these
cases, $(s_1,s_2)_w = s_1c\overline{ s_2}-\overline{s_1}cs_2,\
s_1-as_2b,\ s_1b-as_2$ correspondingly.

For algebras over a field, a composition $(s_1,s_2)_w$ relative to a
trivial lcm $w$ is trivial $mod(s_1,s_2;w)$.

A polynomial $f$ is trivial $mod(S,w)$ if $f$ is a linear
combination of ``normal $S$-words" $(asb),\ s\in S,\ a,b\in N(X)$,
such that $\overline{(asb)} = (a\bar s b)<w$.

The statement $(i)\Rightarrow(ii)$ is the main statement of a
Composition-Diamond lemma since others are much easier to prove. For
two known cases, conformal algebras and dialgebras, (i) and (ii) are
not equivalent.

The case of Lie algebras is an exception since $\bar s\notin N(X)$,
the maximal word of a Lie polynomial $s$ is the maximal word of $s$
as associative polynomial after working out all Lie brackets.

For algebras over a commutative algebra $K$ one needs to deal with a
``double free" algebra $M_{k[Y]}(X)$, a free $M$-algebra over a
polynomial algebra. In this case
$$
lcm(u^Yu^X, v^Yv^X) =  lcm(u^Y, v^Y) lcm(u^X, v^X)
$$
and there are generally infinitely many compositions for given $s_1,
s_2$. It is since we need to use a ``trivial" $lcm(u^X,v^X) =
(u^X)c^X(v^X)$. For algebras over a field, the composition
corresponding to a trivial $(u^X)c^X(v^X)$ is trivial. But for
algebras over a commutative algebra it is not the case if $lcm(u^Y,
v^Y) \neq u^Yv^Y$.

\ \

Recently some new Composition-Diamond lemmas are given: for free
algebra $k\langle Y\rangle\otimes k\langle X\rangle$
\cite{[BCC08]}, for Lie algebras over commutative algebras
\cite{[BCC11]}, for metabelian Lie algebras \cite{[CC12]}, for
semirings \cite{[BCMo13]}, for Rota-Baxter algebras \cite{BCD08},
for L-algebras \cite{[BCH13]}, for Vinberg-Koszul-Gerstenhaber
right-symmetric algebras \cite{BoChLi-GSB-rightsym}, for categories
\cite{BoChLi-CD-category}, for dialgebras \cite{[BCLiu10]}, for
associative algebras with multiple operations
\cite{BoChQiu-CD-Omega}, for associative $n$-conformal algebras
\cite{BoChZhang-n-conf}, for
 associative conformal algebras \cite{BFK04}, for Lie superalgebras \cite{[BokKangLeeMal]},
 for differential algebras \cite{ChChLi-GSB-diff}, for
$\lambda$-differential associative algebras with multiple operators
\cite{ChQiu-Cd-diff}, for commutative algebras with multiple
operators  and free commutative Rota-Baxter algebras \cite{Qiu},
etc.

By using the above and the known Composition-Diamond lemmas,  some
applications are obtained: for embeddings of algebras
\cite{[BCMo10],ChShaoShum}, for free inverse semigroups
\cite{[BCZh09]}, for conformal algebras \cite{BFKK00}, for relative
Gr\"{o}bner-Shirshov bases of algebras and groups \cite{BokutShum},
for extensions of groups and algebras \cite{[Chen08],[Chen09]}, for
some word problems \cite{BokutChainikov08,ChChLuo}, for some Lie
algebras \cite{ChenLiTang}, for partially commutative Lie algebras
\cite{[ChMo-Partialcomm],[Po]}, for braid groups
\cite{[Bok08],[Bok09],[Bo-Ch-Shum],ChMo-GSB-Braid-Artin,[CZhong13]},
for PBW theorems
\cite{[BCC11],BoChLi-GSB-rightsym,[BCLiu10],BoMa99,CCZ,CL,ChenMo-AMS},
etc.

For development of  Gr\"{o}bner-Shirshov bases, one may refer to
surveys: \cite{survey08,BCS,BK03,BK05,[BoMa97]}.

\section{Gr\"{o}bner-Shirshov bases for associative algebras and Lie algebras}

Gr\"{o}bner-Shirshov bases for Lie algebras is established by
Shirshov \cite{[Sh62],[Sh09]} for the free Lie algebras (with
deg-lex order) in 1962 (see also Bokut \cite{b72}). In 1976, Bokut
\cite{[Bok76]} specialized the approach of Shirshov to associative
algebras (see also Bergman \cite{b}). For commutative polynomials,
this is due to Buchberger  \cite{bu65,bu70}.

\subsection{Composition-Diamond lemma for associative algebras}

Let $k\langle X\rangle$ be the free associative algebra over a field
$k$ generated by $X$ and $ (X^{*},<)$ a well-ordered free monoid
generated by $X$, $S\subset k\langle X\rangle$ such that every $s\in
S$ is monic ($s$ is monic if the coefficient of the leading word of
$s$ is 1).

Let us prove $(i)\Rightarrow(ii)$ and define a GS basis.

Let $f=\sum_{i=1}^n\alpha_ia_is_ib_i\in Id(S) $ where each
$\alpha_i\in k, \ a_i,b_i\in X^*, \ s_i\in S$. Let
$w_i=a_i\overline{s_i}b_i, \ w_1=w_2=\cdots=w_l>w_{l+1}\geq\cdots$.

For $l=1$, it is ok.

For $l>1,\ w_1=a_1\overline{s_1}b_1= a_2\overline{s_2}b_2$, common
multiple of $\overline{s_1},\overline{s_2}$, by definition,
$$
w_1=cwd,\ w=``lcm"(\overline{s_1},\overline{s_2}),\
a_is_ib_i=w|_{\overline{s_i}\mapsto s_i}, \ i=1,2,
$$
where $lcm(u,v)\in \{ucv, c\in X^* (a \ trivial \ lcm(u,v));$ $\
u=avb,\ a,b\in X^*\ (an \ inclusion \ lcm(u,v));$ $\ ub=av,\ a,b\in
X^*,\ |ub|<|u|+|v|\ (an \ intersection \ lcm(u,v)\}.$

Then $a_1s_1b_1-a_2s_2b_2=c(w|_{\overline{s_1}\mapsto
s_1}-w|_{\overline{s_2}\mapsto s_2})d=c(s_1,s_2)_wd$. By definition
of GS basis, $(s_1,s_2)_w\equiv0\ \ mod(S,w)$. So,
$a_1s_1b_1-a_2s_2b_2\equiv0\ \ mod(S,w_1)$. We can decrease $l$. By
induction on $l$ and $w_1$, $\bar f=a\bar sb,\ a,b\in X^*,\ s\in S$.

\subsection{Composition-Diamond lemma for Lie algebras over a field}

Let $S\subset{Lie(X)}\subset{k\langle X\rangle}$ be a nonempty set
of monic Lie polynomials, $(X^*,<)$ deg-lex order, $\bar s$ means
the maximal word of $s$ as non-commutative polynomial. Then
compositions are defined as follows
$$
\langle s_1,s_2\rangle_w=[w|_{\overline{s_1}\mapsto
s_1}]_{\overline{s_1}}-[w|_{\overline{s_2}\mapsto
s_2}]_{\overline{s_2}},\ \ w\in ALSW(X)
$$
associative composition with the special Shirshov bracketing, where
$ALSW(X)$ is the set of associative Lyndon-Shirshov words on $X$.

\ \

\noindent{\bf Composition-Diamond lemma for Lie algebras over a
field} (\cite{[BC07],[Sh62],[Sh09]}).  The following statements are
equivalent.
\begin{enumerate}
\item[(i)] $S$ is a Lie GS basis in
$Lie(X)$ (any composition is trivial modulo $(S,w)$).
\item[(ii)] $f\in{Id_{Lie}(S)}\Rightarrow{\bar{f}=a\bar{s}b}$ for
some $s\in{S}$ and $a,b\in{X^*}$.
\item[(iii)]$Irr(S)=\{[u]\in NLSW(X) \ | \  u\neq{a\bar{s}b}, \
s\in{S},\ a,b\in{X^*}\}$ is a linear basis for $Lie(X|S)$.
\end{enumerate}

\section{PBW theorems}

There are 8 PBW (Poincare-Birkhoff-Witt) theorems that are proved by
using GS bases and Composition-Diamond lemmas.

\subsection{Lie algebras--associative algebras (Shirshov)}

Let $L=Lie_k(X|S)$ be an arbitrary Lie algebra with generators $X$
and defining relations $S$ and $ U(L)=k\langle X|S^{(-)}\rangle$ the
universal enveloping associative algebra of $L$. Then
\begin{enumerate}
\item[(i)] $S$ is a Lie GS basis $\Leftrightarrow$ $S$ is an associative GS basis.
\item[(ii)] In this case, a linear basis of $U(L)$ is
\begin{eqnarray*}
u_1u_2\cdots u_t,\ \ u_1\preceq u_2\preceq\cdots\preceq u_t\ \mbox{
(lex-order)},\ \  u_i\in Irr(S)\cap ALSW(X).
\end{eqnarray*}
\end{enumerate}
One uses Shirshov factorization theorem:
$$
u\in X^*,\ \exists! \ u=u_1\cdots u_t,\ \ u_1\preceq \cdots\preceq
u_t,\  u_i\in  ALSW(X).
$$

\subsection{ Lie algebras--pre-Lie algebras (Segal)}

Pre-Lie algebras are defined by an identity $(x,y,z)=(x,z,y)$, where
$(x,y,z)=(xy)z-x(yz)$.

Let $X=\{x_i|i\in I\}$ be a linear basis of the Lie algebra $L$ and
$ [x_i,x_j]=\sum\alpha_{ij}^tx_t=:\{x_i,x_j\} $ the  multiplication
table of the linear basis $X$. Then $L$ has a presentation
$$
L=Lie(x_i,i\in I|[x_i,x_j]=\{x_i,x_j\},\ i,j\in I)=Lie(X|S)
$$
and the universal enveloping pre-Lie algebra of $L$
$$
U_{\mbox{{\small pre-Lie}}}(L)=\mbox{pre-Lie}(X|S^{(-)}).
$$

Then $S^{(-)}$ is a GS basis  of $U_{\mbox{{\small pre-Lie}}}(L)$,
$L\subset U_{\mbox{{\small pre-Lie}}}(L)$ and $Irr(S)$ is a linear
basis of $U_{\mbox{{\small pre-Lie}}}(L)$ by Composition-Diamond
lemma for pre-Lie algebras (Bokut-Chen-Li
\cite{BoChLi-GSB-rightsym}).

\subsection{Leibniz algebras--dialgebras (Aymon,
Grivel)}

Dialgebra: $a\dashv(b\vdash c)=a\dashv b\dashv c,\ (a\dashv b)\vdash
c=a\vdash b\vdash c,\ a\vdash(b\dashv c)=(a\vdash b)\dashv c$ and
$\vdash,\dashv$ associative.

Leibniz identity: $[[a,b],c]=[[a,c],b]+[a,[b,c]]$.

Di-commutator: $[a,b]=a\dashv b-b\vdash a$.
\begin{eqnarray*}
&&L=Lei(x_i,i\in I|[x_i,x_j]=\{x_i,x_j\},\ i,j\in I),\\
&&U_{Dialg}(L)=D(X|S^{(-)}).
\end{eqnarray*}

 A GS basis for $U_{Dialg}(L)$ is given by Bokut-Chen-Liu \cite{[BCLiu10]}
and then a linear basis for $U_{Dialg}(L)$ by Composition-Diamond
lemma for dialgebras which implies  $L\subset U_{Dialg}(L)$.

\subsection{Akivis algebras--non-associative algebras (Shestekov)}

Any nonassociative algebra is an Akivis algebra relative the
commutator $[x,y]=xy-yx$ and the associator $(x,y,z)=(xy)z - x(yz)$.

Akivis identity:
$[[x,y],z]+[[y,z],x]+[[z,x],y]=(x,y,z)+(z,x,y)+(y,z,x)-(x,z,y)-(y,x,z)-(z,y,x)$.
\begin{eqnarray*}
&&A=A(x_i,i\in I|[x_i,x_j]=\{x_i,x_j\},\ (x_i,x_j,x_t)=\{x_i,x_j,x_t\},\ i,j,t\in I),\\
&&U(A)=k\{X|S^{(-)}\},\\
&& S^{(-)}=\{[x_i,x_j]=\{x_i,x_j\},\ (x_i,x_j,x_t)=\{x_i,x_j,x_t\},\
i,j,t\in I\}.
\end{eqnarray*}

A GS basis of $U(A)$ is given by Chen-Li \cite{CL} and then
$A\subset U(A)$.

\subsection{Sabinin algebras--modules (Perez-Izquierdo)}

Let $(V,\langle;\rangle)$ be a Sabinin algebra,
\begin{eqnarray*}
\widetilde{S}(V)&=&T(V)/span\langle xaby-xbay +\sum x_{(1)}\langle
x_{(2)};a,b\rangle y|x,y\in T(V),a,b\in V\rangle\\
&\cong&mod\langle X|I\rangle_{k\langle X\rangle}\ \ \mbox{as }\
k\langle X\rangle\mbox{-modules}
\end{eqnarray*}
the universal enveloping module for $V$, where $I=\{xab-xba+\sum
x_{(1)}\langle x_{(2)};a,b\rangle|x\in X^*,a>b, \ a,b\in X\}$.

Then $I$ is a GS basis (Chen-Chen-Zhong \cite{CCZ}) and then
$V\subset \widetilde{S}(V)$.

\subsection{Dendriform algebras--Rota-Baxter algebras (Chen-Mo, Kolesnikov)
}

Rota-Baxter identity:
\begin{eqnarray*}
P(x)P(y)=P(P(x)y)+P(xP(y)), \forall x,y\in A.
\end{eqnarray*}
Dendriform identities: $ (x\prec y)\prec z=x\prec(y\prec z+y\succ
z),\ (x\succ y)\prec z=x\succ(y\prec z),\ (x\prec y+x\succ y)\succ
z=x\succ(y\succ z). $
\begin{eqnarray*}
D&=&Den(X|x_i\prec x_j=\{x_i\prec x_j\},\  x_i\succ x_j=\{x_i\succ x_j\}, x_i,x_j\in X);\\
U(D)&=&RB(X|x_iP(x_j)=\{x_i\prec x_j\},\  P(x_i)x_j=\{x_i\succ
x_j\}, x_i,x_j\in X).
\end{eqnarray*}

Then $D\subset U(D)$, see Chen-Mo \cite{ChenMo-AMS}.

\subsection{Shirshov's, Cartier's, Cohn's counter
examples to PBW for Lie algebras over commutative algebra}

Shirshov \cite{[Sh53b]} 1953 and Cartier \cite{[Cartier]} 1958 give
counter examples to PBW for Lie algebras over commutative algebra.
Cohn \cite{[Cohn63]} 1963 posts the conjecture:
\begin{eqnarray*}
\mathcal{L}_p&=&Lie_{K}(x_1,x_2,x_3 \ | \ y_3x_3=y_2x_2+y_1x_1),\\
K&=&k[y_1,y_2,y_3|y_i^p=0, i=1,2,3].
\end{eqnarray*}
$\mathcal{L}_p$ can not be embedded into its universal enveloping
associative algebra.

Bokut-Chen-Chen \cite{[BCC11]} establish GS bases theory for Lie
algebras over a commutative algebra. We prove Cohn's conjecture is
correct for $p=2,3,5$.

\subsection{ ``1/2 PBW theorem" (Bokut-Fong-Ke)}

Let $k$ be a field with the characteristic 0 and $C$  a $k$-algebra
with operations $a(n)b,\ n\geq 0$, and $D(a)$. Then $C$ is called a
conformal algebra, if

(1) (the locality condition). For any $a,b\in C$, there exist
$N(a,b)\geq 0$ such that $a(n)b=0$ for $n\geq N(a,b)$.

(2) $D(a(n)b)= D(a)(n)b +a(n)D(b)$ and $D(a)(n)b= -na(n-1)b$.

Let $C$ be a  conformal algebra. Then $C$  is an associative
conformal algebra if
$$
(As)\ \ (a(n)b)(m)c= \sum_{s\geq 0}(-1)^sC{_n}^s a(n-s)(b(m+s)c);
$$
$C$ is a Lie conformal algebra if
\begin{eqnarray*}
&&(Lie)\ \  \mbox{(anti-commutativity) }\ (a[n]b)=\{b[n]a\}, \\
&& \ \ \ \ \mbox{ where }\ \{b[n]a\} = \sum_{s\geq
0}(-1)^{n+s}D^{(s)}(b[n+s]a),\
D^{(s)}=1/s!D^s,\\
&& \mbox{(Jacoby) }\ ((a[n]b)[m]c) = \sum_{s\geq
0}(-1)^{s}C_{n}^s((a[n-s](b[m+s]c))- (b[m+s](a[n-s]c))).
\end{eqnarray*}

 Gr\"{o}bner-Shirshov bases for conformal associative algebras ($n$-conformal associative
 algebras) are established in Bokut-Fong-Ke \cite{BFK04} and
Bokut-Chen-Zhang \cite{BoChZhang-n-conf}.

There is a 1/2 PBW-theorem for Lie conformal algebras: Let
$$
L=LieCon(\{a_i,\ i\in I\}, N| a_i[n]a_j =\sum\alpha_{ij}^k a_k,\
i\geq j,\ n<N)
$$
be a Lie conformal algebra with the linear basis $\{a_i,\ i\in I\}$
over $k[D]$ and the locality $N$. Let
$$
U(L) = AsCon(\{a_i,\ i\in I\}, N| S=\{s_{ij}^n= a_i(n)a_j -
\{a_j(n)a_i\} - a_i[n]a_j,\ i\geq j,\ n<N\})
$$
be the universal enveloping associative conformal algebra for $L$.

Then any composition $(s_{ij}^n, s_{jk}^m)_w=s_{ij}^n(m)a_k -
a_i(n)s_{jk}^m,\ i>j>k,\ n,m<N$ is trivial $mod(S,w),\
w=a_i(n)a_j(m)a_k$.

\ \

There is also a 1/2 PBW theorem between $n$-conformal Lie algebras
and $n$-conformal associative algebras, see \cite{BoChZhang-n-conf}.

\section{Linear bases of free universal algebras}

--Bases of free Lie algebras

M. Hall and  A.I. Shirshov use constructions and check axioms.

One may use anti-commutative Gr\"obner-Shirshov bases for a free Lie
algebra.

Hall basis (Bokut-Chen-Li \cite{[BCLi09]}): $Lie(X)=AC(X|S_1)$,
$S_1$ is a anti-commutative GS basis, $Irr(S_1)=$Hall basis in $X$.

Lyndon-Shirshov basis (Bokut-Chen-Li \cite{[BCLi13]}):
$Lie(X)=AC(X|S_2)$, $S_2$ is a anti-commutative GS basis,
$Irr(S_2)=$Lyndon-Shirshov basis in $X$.

\ \

--Loday basis of a free dialgebra

$D(X)=L(X|S)$, $L$-identity: $(a\vdash b)\dashv c=a\vdash (b\dashv
c)$, $S$ an $L$-GS basis (GS basis as $L$-algebra),  $Irr(S)=$Loday
basis in $X$ (Bokut-Chen-Huang \cite{[BCH13]}).

\ \

 --Bases of a free dentriform algebra

$Den(X)=L(X|S)$, $Irr(S)$=a linear basis of $Den(X)$
(Bokut-Chen-Huang \cite{[BCH13]}).

\ \

--Bases of a free Rota-Baxter algebra (Rota, Cartier)

Via GS method for $\Omega$-algebras (Bokut-Chen-Qiu
\cite{BoChQiu-CD-Omega}).

\ \

--Free inverse semigroup (Polyakova-Schain)

An associative GS basis is given by (Bokut-Chen-Zhao
\cite{[BCZh09]}), $Irr(S)$ is a normal form of free inverse
semigroup.

\ \

--Free idempotent semigroup (Chen-Yang \cite{Chen-Yang}).

\section{Normal forms for groups and semigroups}

--Braid groups

in Artin-Burau generators (Bokut-Chanikov-Shum \cite{[Bo-Ch-Shum]});

in Artin-Garside generators (Bokut \cite{[Bok08]});

in Birman-Ko-Lee generators (Bokut \cite{[Bok09]});

in Adyan-Thurston generators (Chen-Zhong \cite{[CZhong13]}).

\ \

--Chinese monoid (Chen-Qiu \cite{Chen-Qiu})

\ \

--Plactic monoid (Bokut-Chen-Chen-Li \cite{BoChLiJing-plactic}).

\ \

--HNN extension

Britton Lemma and Lyndon-Schupp normal form lemma for  an
HNN-extension of a group was proved using an associative
Composition-Diamond lemma relative to a ``$S$-partially" monomial
order of words (Chen-Zhong \cite{[ChenZhong08]}).

\ \

--One-relator groups

In (Chen-Zhong \cite{[CZhong11]}), some one-relator groups were
studying by means of groups with the standard normal forms (the
standard  GS bases) in the sense (Bokut \cite{[Bok66],[Bok67]}). It
is known that any one-relator group can be effectively embedded into
2-generator one-relator group $G=gp\langle x,y| x^{i_1}y^{j_1}\dots
x^{i_k}y^{j_k}\rangle,\ k\geq 1$ is the depth. It is proved that  a
group $G$ of depth $\leq 3$ is computably embeddable into  a
 Magnus-Moldavanskii tower of HNN-extensions with the standard normal form
of elements. There are quite a lot of examples that support an old
conjecture that the result is valid for any depth.

\section{Schreier extensions of groups and algebras}

In (Chen \cite{[Chen08]}), it is dealing with a Schreier extension
$$
1\rightarrow A\rightarrow C\rightarrow B\rightarrow 1
$$
of a group $A$ by $B$. M. Hall \cite{[H]} mentioned that for any $B$
it is difficult to find an analogous conditions. Actually this
problem was solved in \cite{[Chen08]} using the GS bases technique.
As applications there were given above conditions for cyclic and
free abelian cases, as well as for the case of HNN-extensions.

The same kind of result was established for Schreier extensions of
associative algebras (Chen  \cite{[Chen09]}).

\section{Embedding of algebras}

In Bokut-Chen-Mo \cite{[BCMo10]}, we were dealing with the problem
of embedding of countably generated associative and Lie algebras,
groups, semigroups, $\Omega$-algebras into (simple) 2-generated
ones. We proved some known results (of Higman-Neuman-Neuman, Evance,
Malcev, Shirshov) and some new ones using GS bases technique. For
example

\ \

\noindent{\bf Theorem 1}. Every countable Lie algebra is embeddable
into simple 2-generated Lie algebra.

\ \

\noindent{\bf Theorem 2}. Every countable differential algebra is
embeddable into a simple 2-generated differential algebra.

\ \

G. Bergman (Private communication, 2013 \cite{[Ber]}) gave an idea
how to avoid the restriction on cardinality of  the ground field.
Now Qiuhui Mo proved that the Bergman's idea works.

\ \

\ \

\ \

\end{document}